\documentclass[12pt]{amsart}

\newcommand {\supplus}{\mathop{{\supset}\llap{\raise 
0.5pt\hbox{\normalfont\small+}\hskip 0.5pt}}} 

\newcommand {\subplus}{\mathop{{\subset}\llap{\raise 
0.5pt\hbox{\normalfont\small+}\hskip 0.5pt}}}  

\newcommand {\cal} {\mathcal}

\def \opname#1#2%
  {\expandafter\newcommand \csname #1\endcsname {{\mathop{#2}\nolimits}}}

\newcommand{\rmname}[1]
  {\expandafter\newcommand \csname #1\endcsname {{\operatorname{#1}}}}

\newcommand{\rmnameii}[2]
  {\expandafter\newcommand \csname #1\endcsname {{\operatorname{#2}}}}

\rmname{act}
\rmname{Ad}
\rmname{Add}
\rmname{ad}
\rmname{Alt}
\rmname{alt}
\rmname{Ann}
\rmname{antidiag}
\rmname{Ber}
\rmname{ber}
\rmname{Br}
\rmname{card}
\rmname{ch}
\rmname{Char}
\rmname{cem}
\rmname{cj}
\rmname{Cliff}
\rmname{cntr}
\rmname{codim}
\rmname{coind}
\rmname{const}
\rmname{col}
\rmname{cork}
\rmname{cpr}
\rmname{diag}
\rmnameii{Div}{div}
\rmname{Def}
\rmname{Der}
\rmname{Dim}
\rmname{End}
\rmname{Even}
\rmname{Ext}
\rmname{gr}
\rmname{Hom}
\rmname{HT}
\rmnameii{Ht}{ht}
\rmname{hwt}
\rmname{Id}
\rmname{id}
\rmname{ind}
\rmname{Ind}
\rmname{Inf}
\rmname{irr}
\rmname{Le}
\rmname{Lie}
\rmname{lwt}
\rmname{mult}
\rmname{Mor}
\rmname{nm}
\rmname{Ob}
\rmname{Odd}
\rmname{Osc}
\rmname{per}
\rmname{Pic}
\rmname{pr}
\rmname{pro}
\rmname{Prime}
\rmname{Proj}
\rmname{prt}
\rmname{pt}
\rmname{Q}
\rmname{qet}
\rmname{qtr}
\rmname{rd}
\rmname{rk}
\rmname{row}
\rmname{Res}
\rmname{salt}
\rmname{Sch}
\rmname{SBr}
\rmname{scalar}
\rmname{Ser}
\rmname{sign}
\rmname{Smbl}
\rmname{spin}
\rmname{ssym}
\rmname{str}
\rmname{st}
\rmname{sgn}
\rmname{sq}
\rmname{symm}
\rmname{supp}
\rmname{Supp}
\rmname{St}
\rmname{Spec}
\rmname{Spm}
\rmname{tr}
\rmname{vpt}
\rmname{weyl}
\rmname{Weyl}
\rmname{Witt}

\opname{vvol}  {{v\hspace{-0.1ex}o\hspace{-0.02ex}l\/}}
\opname{pnt}  {\text{\normalfont pt}}
\opname{Span} {{Span}}
\opname{slim} {\overline{\lim}}
\opname{Vol}  {{V\hspace{-0.55ex}o\hspace{-0.02ex}l\/}}
\opname{Par}  {{P\hspace{-0.3ex}a\hspace{-0.05ex}r\/}}

\rmname{Mat}
\rmname{Bil}
\rmname{Diff}
\rmname{Ker}
\rmname{Herm}
\rmname{Coker}
\rmname{Conn}
\rmname{Covect}
\rmname{Vect}
\rmname{Int}

\rmnameii {IM} {Im}
\rmnameii {RE} {Re}

\opname{Aut} {{A\hspace{-0.2ex}u\hspace{-0.1ex}t\/}}
\opname{GL} {{G\hspace{-0.3ex}L}}
\opname{SL} {{S\hspace{-0.3ex}L}}
\opname{Exp} {{E\hspace{-0.2ex}x\hspace{-0.1ex}p\/}}
\opname{GQ} {{G\hspace{-0.2ex}Q}}
\opname{OSp} {{O\hspace{-0.25ex}S\hspace{-0.15ex}p\/}}
\opname{Out} {{O\hspace{-0.25ex}u\hspace{-0.15ex}t\/}}
\opname{Spp} {{S\hspace{-0.2ex}p\/}}
\opname{SpO} {{S\hspace{-0.2ex}p\hspace{-0.02ex}O\/}}
\opname{Pe} {{P\hspace{-0.25ex}e\/}}
\opname{SPe} {{S\hspace{-0.25ex}P\hspace{-0.25ex}e\/}}
\opname{Spin} {{S\hspace{-0.25ex}p\hspace{-0.05ex}i\hspace{-0.1ex}n\/}}
\opname{Iso} {{I\hspace{-0.25ex}s\hspace{-0.1ex}o\/}}
\opname{SSPe} {{S\hspace{-0.25ex}S\hspace{-0.15ex}P\hspace{-0.25ex}e\/}}
\opname{PeU} {{P\hspace{-0.25ex}e\hspace{-0.1ex}U\/}}
\opname{QU} {{Q\hspace{-0.15ex}U\/}}
\opname{U} {{U\/}}

\opname{cGQ} {{\cal G \hspace{-0.2em} Q \/}}
\opname{cSL} {{\cal S \hspace{-0.2em} L \/}}
\opname{cGL} {{\cal G \hspace{-0.2em} L \/}}
\opname{cOSp} {{\cal O \hspace{-0.2em} S \hspace{-0.3em} \it p\/}}
\opname{cPe} {{\cal P \hspace{-1.5pt} \it e\/}}
\opname{cVect} {{\cal V \hspace{-1.5pt} \it e\hspace{-0.1ex}c\hspace{-0.1ex}t\/}}
\opname{cVol} {{\cal V \hspace{-1.5pt} \it o\hspace{-0.1ex}l\/}}
\opname{cAut} {{\cal A \hspace{-0.2em} \it u\hspace{-0.1em}t\/}}
\opname{cCovect} {{\cal C \hspace{-1.5pt}
     \it o\hspace{-0.1ex}v\hspace{-0.1ex}e\hspace{-0.1ex}c\hspace{-0.1ex}t\/}}
\opname{CW} {{C\hspace{-0.15ex}W}}

\newcommand {\bcdot}   {\mathbin{\hbox{\raise.4ex\hbox{\bf.}}}} 

\newcommand {\secno} {}

\newtheorem{Theorem}{\secno Theorem}

\newenvironment {th*}[1]
    {\gdef\thname{#1} \begin{thn}}%
    {\end{thn}}
\newtheorem{thn}[Theorem] {\thname}

\theoremstyle{definition}

\newenvironment {ex*}[1]
    {\gdef\thname{#1} \begin{exn}}%
    {\end{exn}}
\newtheorem{exn}[Theorem]{\thname}

\theoremstyle{remark}

\newenvironment {rem*}[1]
    {\gdef\thname{#1} \begin{remn}}%
    {\end{remn}}
\newtheorem{remn}[Theorem]{\thname}

\begin{document}

\title{The Howe duality and Lie superalgebras}

\author{Dimitry Leites${}^1$, Irina 
Shchepochkina${}^2$} 

\address{${}^1$Deptartment of Mathematics, University of Stockholm,
Roslagsv.  101, Kr\"aftriket hus 6, S-106 91, Stockholm, Sweden;
mleites@matematik.su.se; ${}^2$The Independent University of Moscow,
Bolshoj Vlasievsky per, dom 11, RU-121 002 Moscow, Russia;
Ira@Pa\-ra\-mo\-no\-va.mccme.ru}



\begin{abstract} Howe's duality is considered from a unifying point of
view based on Lie superalgebras.  New examples are offered.
In particualr, we construct several simplest spinor-oscillator
representations and compute their highest weights for the \lq\lq
stringy" Lie superalgebras (i.e., Lie superalgebras of complex vector
fields (or their nontrivial central extensions) on the supercircle
$S^{1|n}$ and its two-sheeted cover associated with the M\"obius
bundle).
\end{abstract}

\thanks{D.L. is thankful to P. Deligne
whose question prompted \cite{del/lei} and this paper, to B.~Feigin, 
E.~Poletaeva, V.~Serganova and Xuan Peiqi for help; we gratefully
acknowledge financial support of an NFR grant and RFBR grant
99-01-00245, respectively.}

\maketitle

In our two lectures we briefly review, on the most elementary level,
several results and problems unified by ``Howe's duality''.  Details
will be given elsewhere.  The ground field in the lectures is
${\mathbb C}$.

\section*{\S 1. Introduction}
 
In his famous preprint \cite{how1} R.~Howe gave an inspiring
explanation of what can be \lq\lq dug out" from H.~Weyl's \lq\lq
wonderful and terrible" book \cite{wey}, at least as far as invariant
theory is concerned, from a certain unifying viewpoint.  According to
Howe, much is based on a remarkable correspondence between certain
irreducible representations of Lie subalgebras $\Gamma$ and
$\Gamma'$ of the Lie algebra $\mathfrak{o}(V)$ or $\mathfrak{sp}(V)$
provided $\Gamma$ and $\Gamma '$ are each other's \lq\lq commutants",
i.e., centralizers.  This correspondence is known ever since as {\it
Howe's correspondence} or {\it Howe's duality}.\index{Howe's
duality}\index{Howe's correspondence} In \cite{how1} and subsequent
papers Howe gave several examples of such a correspondence previously
known, mostly, inadvertently.  Let us remind some of 
them (omitting important Jacquet-Langlands-Shimizu 
correspondence, S.~Gelbart's contributions, etc.) :

1) decomposition of $\mathfrak{o} (V)$-module $S^{\bcdot}(V)$ into
spherical harmonics;

2) {\it Lefschetz decomposition} \index{Lefschetz decomposition} of
$\mathfrak{sp}(V)$-module $\Lambda ^{\bcdot}(V)$ into primitive forms
(sometimes this is called Hodge--L\'epage decomposition);

3) a striking resemblance between {\it spinor
representation}\index{representation, spinor} of $\mathfrak{o} (n)$
and {\it oscillator} ({\it Shale--Segal--Weil--metaplectic}--...  )
\index{representation, oscillator}\index{representation,
Shale--Segal--Weil--metaplectic--...} representation of $\mathfrak{sp}
(2n)$.

As an aside Howe gives the \lq\lq shortest possible" proof of the {\it
Poincar\'e lemma}.\index{Poincar\'e lemma}\index{lemma, Poincar\'e}
(Recall that this lemma states that in any sufficiently small open
star-shaped neighborhood of any point on any manifold any closed
differential form is exact.)  In this proof, Lie superalgebras, that
lingered somewhere in the background in the previous discussion but
were treated rather as a nuisance than help, are instrumental to reach
the goal.  This example shows also that the requirement of reductivity
of $\Gamma$ and $\Gamma'$ to form a ``dual pair'' is extra.  Elsewhere
we will investigate what are the actual minimal restrictions on
$\Gamma$ and $\Gamma'$ needed to reach one of the other problems
usually solved by means of Howe duality: decompose the symmetric or
exterior algebra of a module over $\Gamma\oplus \Gamma'$.
Howe's manuscript was written at the time when supersymmetry theory
was being conceived.  By the time \cite{how1} was typed, the
definition of what is nowadays called {\it
superschemes}\index{superscheme} (\cite{lei0}) was not yet rewritten
in terms to match physical papers (language of points was needed; now we can
recommend \cite{del}) nor translated into English and, therefore, was
unknown; the classification of simple finite dimensional Lie
superalgebras over ${\mathbb C}$ had just been announced.  This was,
perhaps, the reason for a cautious tone with which Howe used Lie
superalgebras, although he made transparent how important they might
be for a lucid presentation of his ideas and explicitly stated so.

Since \cite{how2}, the published version of \cite{how1}, though put
aside to stew for 12 years, underwent only censorial changes, we
believe it is of interest to explore what do we gain by using Lie
superalgebras from the very beginning (an elaboration of other aspects
of this idea \cite{del/lei} are not published yet).  Here we briefly
elucidate some of Howe's results and notions and give several new
examples of Howe's dual pairs.  In the lectures we will review the known
examples 1) -- 3) mentioned above but consider them in an appropriate
``super'' setting, and add to them:

4) a refinement of the Lefschetz decomposition --- J.~Bernstein's
decomposition (\cite{ber}) of the space $\Omega_{\hbar}^{\bcdot}$ of
\lq\lq twisted" differential forms on a symplectic manifold with
values in a line bundle with connection whose curvature form differs
by a factor $\hbar$ from the canonical symplectic form;

5) a decomposition of the space of differential forms on a hyper-K\"
ahlerian manifold similar to the Lefschetz one (\cite{ver}) but with
$\mathfrak{sp} (4)$ instead of $\mathfrak{sp} (2)=\mathfrak{sl} (2)$
and its refinement associated with the $\mathfrak{osp} (1|4)$.

6) Apart from general clarification of the scenery and new examples
even in the old setting, i.e., on manifolds, the superalgebras
introduced {\it ab ovo} make it manifest that there are at least two
types of Howe's correspondence: the conventional one and several
\lq\lq ghost" ones associated with quantization of the antibracket
\cite{lei/shc}.

7) Obviously, if $\Gamma\oplus \Gamma'$ is a maximal subalgebra of
$\mathfrak{osp}$, then $(\Gamma, \Gamma')$ is an example of Howe dual
pair.  Section 6 gives some further examples, partly borrowed from
\cite{shc}, where more examples can be found.

We consider here only finite dimensional Lie superalgebras with the
invariant theory in view.  In another lecture (\S\S3,4) we consider
spinor-oscillator representations in more detail.  In these elementary
talks we do not touch other interesting applications such as Capelli
identities (\cite{kos},\cite{mol}), or prime characteristic
(\cite{ryb}).  Of dozens of papers with examples of Howe's duality in
infinite dimensional cases and still other examples, we draw attention
of the reader to the following selected ones: \cite{fei/sem/tip}, and
various instances of bose-fermi correspondence, cf.  \cite{Fre} and
\cite{Kac}.  Observe also that the Howe duality often manifests itself
for $q$-deformed algebras, e.g., in Klimyk's talk at our conference,
or \cite{Din/Fre}.  To treat this $q$-Howe duality in a similar way,
we first have to explicitly $q$-quantize Poisson superalgebras
$\mathfrak{po}(2n|m)$ (for $mn=0$ this is straightforward replacement
of (super)commutators from \cite{Lei/Pol} with $q$-(super)commutators.

\section*{\S 2. The Poisson superalgebra
$\mathfrak{g}=\mathfrak{po}(2n|m)$}

{\bf 2.1.  Certain ${\mathbb Z}$-gradings of $\mathfrak{g}$}.  Recall
that $\mathfrak{g}$ is the Lie superalgebra whose superspace is
${\mathbb C}[q, p, \Theta]$ and the bracket is the {\it Poisson
bracket} $\{\cdot , \cdot\}_{P.b.}$ is given by the formula
$$
\renewcommand{\arraystretch}{1.4}
\begin{array}{l}
    \{f, g\}_{P.b.}=\sum\limits_{i\leq n}\ \bigg(\frac{\partial 
f}{\partial p_i}\ 
\frac{\partial g}{\partial q_i}-\ \frac{\partial f}{\partial q_i}\ 
\frac{\partial g}{\partial p_i}\bigg)-\\
(-1)^{p(f)}\sum\limits_{j\leq m}\ 
\frac{\partial f}{\partial \theta_j}\ 
\frac{\partial g}{\partial \theta_j}\; {\rm  for }\; f, g\in {\mathbb 
C} [p, q, \Theta].
\end{array}\eqno{(2.1)}
$$
Sometimes it is more convenient to redenote the $\Theta$'s and set
$$
\renewcommand{\arraystretch}{1.4}
\begin{array}{l}
    \xi_j=\frac{1}{\sqrt{2}}(\Theta_{j}-i\Theta_{r+j});\quad
\eta_j=\frac{1}{ \sqrt{2}}(\Theta_{j}+i\Theta_{r+j})\\ 
{\rm  for}\;
j\leq r= [m/2]\; ({\rm here}\; i^2=-1), \quad \theta =\Theta_{2r+1}
\end{array}
$$ 
and accordingly modify the bracket (if $m=2r$, there is no term with
$\theta$):
$$
\renewcommand{\arraystretch}{1.4}
\begin{array}{l}
	\{f, g\}_{P.b.}=\sum\limits_{i\leq n}\ \bigg(\frac{\partial f}{\partial p_i}\ 
\frac{\partial g}{\partial q_i}-\ \frac{\partial f}{\partial q_i}\ 
\frac{\partial g}{\partial p_i}\bigg)-\\
(-1)^{p(f)}\bigg[\sum\limits_{j\leq m}(\frac{\partial f}{\partial\xi_j}
\frac{\partial g}{\partial \eta_j}+\frac{\partial f}{\partial \eta_j}\ 
\frac{\partial g}{\partial \xi_j})+\frac{\partial f}{\partial \theta}\ 
\frac{\partial g}{\partial \theta}\bigg].
\end{array}
$$
Setting $\deg_{Lie} f=\deg f-2$ for any monomial $f\in {\mathbb C} [p,
q, \Theta]$, where $\deg p_{i}=\deg q_{i}=\deg \Theta _{j}=1$ for all
$i, j$, we obtain the {\it standard} ${\mathbb Z}$-grading of
$\mathfrak{g}$:

\begin{tabular}{|c|c|c|c|c|c|}
\hline
degree of $f$&$-2$&$-1$&$0$&$1$&$\dots$\\
\hline
$f$&$1$&$p, \, q, \, \theta $&$f: \deg f=2$&$f: \deg
f=3$&$\dots$\\
\hline
\end{tabular}

\noindent Clearly, $\mathfrak{g}=\mathop{\oplus}\limits_{i\geq
-2}\mathfrak{g}_i$ with $\mathfrak{g}_0\simeq\mathfrak{osp}(m|2n)$. 
Consider now another, ``{\it rough}'', grading of $\mathfrak{g}$.  To
this end, introduce: $Q=(q, \xi)$, $P=(p, \eta)$ and set
$$
\deg Q_i=0,\; \deg \theta=1,\; \deg P_i=\left\{\begin{array}{ll}1& {\rm if} \; 
m=2k\\
2& {\rm if}\; m=2k+1.\end{array}\right.\eqno{(*)}
$$
{\bf Remark}.  Physicists prefer to use half-integer values of $\deg$
for $m=2k+1$ by setting $\deg \theta=\frac12$ and $\deg P_i=1$ at all
times.

The above grading $(*)$ of the polynomial algebra induces the
following {\it rough grading} of the Lie superalgebra $\mathfrak{g}$. 
For $m=2k$ just delete the columns of odd degrees and delete the 
degrees by 2:

\noindent $m=2k+1$: \begin{tabular}{|c|c|c|c|c|c|c|}
\hline
degree&$\dots$&$2$&$1$&$0$&$-1$&$-2$\\
\hline
elements&$\dots$&${\mathbb C} [Q]P^2$&${\mathbb C} [Q]P\theta$& ${\mathbb C} [Q]P$&${\mathbb C} 
[Q]\theta$&${\mathbb C} [Q]$\\
\hline
\end{tabular}

{\bf 2.2.  Quantization}.  We call the nontrivial deformation
${\mathcal Q}$ of the Lie superalgebra $\mathfrak{po}(2n|m)$ {\it
quantization}\index{quantization} (for details see \cite{lei/shc}). 
There are many ways to quantize $\mathfrak{g}$, but all of them are
equivalent.  Recall that we only consider $\mathfrak{g}$ whose
elements are represented by polynomials; for functions of other types
(say, Laurent polynomials) the uniqueness of quantization may be
violated.

Consider the following quantization, so-called $QP$-quantization, 
given on linear terms by the formulas:
$$
{\mathcal Q}: Q\mapsto \hat Q, \quad P\mapsto 
\hbar\frac{\partial}{\partial Q},\eqno{(*)}
$$
where $\hat Q$ is the operator of left multiplication by $Q$; an
arbitrary monomial should be first rearranged so that the $Q$'s
stand first (normal form) and then apply $(*)$ term-wise.

The deformed Lie superalgebra ${\mathcal Q}(\mathfrak{po}(2n|2k))$ is
the Lie superalgebra of differential operators with polynomial
coefficients on ${\mathbb R}^{n|k}$.  Actually, it is an analog of
$\mathfrak{gl}(V)$.  This is most clearly seen for $n=0$. 
Indeed,
$$
{\mathcal
Q}(\mathfrak{po}(0|2k))=\mathfrak{gl}(\Lambda^{\bcdot}(\xi))=
\mathfrak{gl}(2^{k-1}|2^{k-1}).
$$
In general, for $n\neq 0$, we have
$$
{\mathcal Q}(\mathfrak{po}(2n|2k))={\rm ``}\mathfrak{gl}{\rm "}
({\mathcal F}(Q))= \mathfrak{diff}({\mathbb R}^{n|k}).
$$
For $m=2k-1$ we consider $\mathfrak{po}(0|2k-1)$ as a subalgebra of
$\mathfrak{po}(0|2k)$; the quantization sends $\mathfrak{po}(0|2k-1)$
into $\mathfrak{q}(2^{k-1})$.  For $n\neq 0$ the image of ${\mathcal
Q}$ is an infinite dimensional analog of $\mathfrak{q}$, indeed (for
$J=i(\theta+\frac{\partial}{\partial \theta})$ with $i^2=-1$):
$$
{\mathcal Q}(\mathfrak{po}(2n|2k-1))=\mathfrak{qdiff}({\mathbb
R}^{n|k})=
\{D\in\mathfrak{diff}({\mathbb R}^{n|k}): [d, J]=0 \}.
$$

{\bf 2.3.  Fock spaces and spinor-oscillator representations}.  The
Lie superalgebras $\mathfrak{diff}({\mathbb R}^{n|k})$ and
$\mathfrak{qdiff}({\mathbb R}^{n|k})$ have indescribably many
irreducible representations even for $n=0$.  But one of the
representations, the identity one, in the superspace of functions on
${\mathbb R}^{n|k}$, is the ``smallest'' one.  Moreover, if we
consider the superspace of $\mathfrak{diff}({\mathbb R}^{n|k})$ or
$\mathfrak{qdiff}({\mathbb R}^{n|k})$ as the {\it associative}
superalgebra (denoted ${\rm Diff}({\mathbb R}^{n|k})$ or ${\rm
QDiff}({\mathbb R}^{n|k})$), this associative superalgebra has only {\it one}
irreducible representation --- the same identity one.  This
representation is called {\it the Fock space}.\index{Fock space}

As is known, the Lie superalgebras $\mathfrak{osp}(m|2n)$ are rigid
for $(m|2n)\neq (4|2)$.  Therefore, the through map
$$
\mathfrak{h}\longrightarrow  \mathfrak{g}_{0}=\mathfrak{osp}(m|2n)\subset
\mathfrak{g}=\mathfrak{po}(2n|m)\stackrel{{\mathcal
Q}}{\longrightarrow }\mathfrak{diff}({\mathbb  R}^{n|k})
$$ 
sends any subsuperalgebra $\mathfrak{h}$ of $\mathfrak{osp}(m|2n)$
(for $(m|2n)\neq (4|2)$) into its isomorphic image.  (One can also
embed $\mathfrak{h}$ into $\mathfrak{diff}({\mathbb R}^{n|k})$
directly.)  The irreducible subspace of the Fock space which contains
the constants is called the {\it spinor-oscillator
representation}\index{spinor-oscillator
representation}\index{representation, spinor-oscillator} of
$\mathfrak{h}$.  In particular cases, for $m=0$ or $n=0$ this subspace
turns into the usual {\it spinor}\index{spinor
representation}\index{representation, spinor} or {\it oscillator
representation}\index{oscillator representation}\index{representation,
oscillator}, respectively.  We have just given a unified description
of them.  (A more detailed description follows.)

{\bf 2.4.  Primitive alias harmonic elements}. The elements of
$\mathfrak{osp}(m|2n)$ (or its subalgebra $\mathfrak{h}$) act in the
space of the spinor-oscillator representation by inhomogeneous
differential operators of order $\leq 2$ (order is just the filtration
associated with the ``rough'' grading):

\begin{tabular}{|c|c|c|c|}
\hline
$m=2k$: &&&\\
\hline
${\rm  degree}$&$-1$&$0$&$1$\\
\hline
${\rm  elements}$&$\hat P^2$&$\hat P\hat Q$&$\hat Q^2$\\
\hline
\end{tabular}
\begin{tabular}{|c|c|c|c|c|c|}
\hline
$m=2k+1$: &&&&&\\
\hline
${\rm  degree}$&$-2$&$-1$&$0$&$1$&$2$\\
\hline
${\rm  elements}$&$\hat P^2$&$\hat P\hat \theta$&$\hat P\hat Q$
&$\hat Q\hat \theta$&$\hat Q^2$\\
\hline
\end{tabular}

The elements from $({\mathbb C} [Q])^{\hat P^{2}}$ for $m=2k$ or
$({\mathbb C} [Q, \theta])^{\hat P\hat \theta}$ for $m=2k+1$ are
called {\it primitive} or {\it harmonic} ones.  
More generally, let $\mathfrak{h}\subset\mathfrak{osp} (m|2n)$ be a
${\mathbb Z}$-graded Lie superalgebra embedded consistently with the
rough grading of $\mathfrak{osp} (m|2n)$.  Then the elements from
$({\mathbb C} [Q])^{\mathfrak{h}_{-1}}$ for $m=2k$ or $({\mathbb C}
[Q, \theta])^{\mathfrak{h}_{-1}}$ for $m=2k+1$ will be called {\it
$\mathfrak{h}$-primitive}\index{primitive element}\index{element,
primitive} or {\it $\mathfrak{h}$-harmonic}\index{harmonic
element}\index{element,harmonic}.

{\bf 2.4.1.  Nonstandard ${\mathbb Z}$-gradings of
$\mathfrak{osp}(m|2n)$}.  It is well known that one simple Lie
superalgebra can have several nonequivalent Cartan matrices and
systems of Chevalley generators, cf.  \cite{Gro/Lei2}.  Accordingly,
the corresponding divisions into {\it positive} and {\it negative}
root vectors are distinct.  The following problem arises: {\sl How the
passage to nonstandard gradings affects the highest weight of the
spinor-oscillator representation defined in sec. 3?} (Cf. 
\cite{nis/hay})

{\bf 2.5.  Examples of dual pairs}.  Two subalgebras $\Gamma, \Gamma'$
of $\mathfrak{g}_{0}=\mathfrak{osp}(m|2n)$ will be called a {\it dual
pair}\index{dual pair}\index{pair, dual} if one of them is the
centralizer of the other in $\mathfrak{g}_{0}$.

If $\Gamma\oplus\Gamma'$ is a maximal subalgebra in
$\mathfrak{g}_{0}$, then, clearly, $\Gamma, \Gamma'$ is a dual pair. 
A generalization: consider a pair of mutual centralizers $\Gamma,
\Gamma'$ in $\mathfrak{gl}(V)$ and embed $\mathfrak{gl}(V)$ into
$\mathfrak{osp}(V\oplus V^*)$.  Then $\Gamma, \Gamma'$ is a dual pair
(in $\mathfrak{osp}(V\oplus V^*)$).  For a number of such examples see
\cite{Shc1}.  Let us consider several of these examples in detail.

{\bf 2.5.1}.  $\Gamma=\mathfrak{sp}(2n)=\mathfrak{sp}(W)$ and
$\Gamma'=\mathfrak{sp}(2)=\mathfrak{sl}(2)=\mathfrak{sp}(V\oplus
V^*)$.  Clearly, $\mathfrak{h}= \Gamma\oplus \Gamma'$ is a maximal
subalgebra in $\mathfrak{o}(W\otimes(V\oplus V^*))$.  The Fock space
is just $\Lambda^{\bcdot}(W)$.

The following classical theorem and its analog 5.2 illustrate the importance
of the above notions and constructions.

{\bf Theorem}. {\it The $\Gamma'$-primitive elements of
$\Lambda^{\bcdot}(W)$ of each degree $i$ constitute an irreducible
$\Gamma$-module $P^i_{\mathfrak{sp}}$, $0\leq i\leq n$.}

This action of $\Gamma'$ in the superspace of differential forms on
any symplectic manifold is well known: $\Gamma'$ is generated (as a
Lie algebra) by operators $X_{+}$ of left multiplication by the
symplectic form $\omega$ and $X_{-}$, application of the bivector dual
to $\omega$.

{\bf 2.5.2}. $\Gamma=\mathfrak{o}(2n)=\mathfrak{o}(W)$ and
$\Gamma'=\mathfrak{sp}(2)=\mathfrak{sl}(2)=\mathfrak{sp}(V\oplus
V^*)$.  Clearly, $\mathfrak{h}= \Gamma\oplus \Gamma'$ is a maximal
subalgebra in $\mathfrak{sp}(W\otimes(V\oplus V^*))$.  The Fock space
is just $S^{\bcdot}(W)$.

{\bf Theorem}. {\it The $\Gamma'$-primitive elements of
$S^{\bcdot}(W)$ of each degree $i$ constitute an irreducible
$\Gamma$-module $P^i_{\mathfrak{o}}$, $i=0, 1, \ldots$. }

This action of $\Gamma'$ in the space of polynomial functions on any
Riemann manifold is also well known: $\Gamma'$ is generated (as a Lie
algebra) by operators $X_{+}$ of left multiplication by the quadratic
polynomial representing the metric $g$ and $X_{-}$ is the
corresponding Laplace operator.

Clearly, a mixture of Examples 2.5.1 and 2.5.2 corresponding to
symmetric or skew-symmetric forms on a supermanifold is also possible:
{\sl the space of $\Gamma'$-primitive elements of $S^{\bcdot}(W)$ of
each degree $i$ is an irreducible $\Gamma$-module}, cf. 
\cite{nis/hay} and Sergeev's papers \cite{ser1}, \cite{ser2}.

In \cite{how1}, \cite{how2} the dual pairs had to satisfy one more
condition: the through action of both $\Gamma$ and $\Gamma'$ on the
identity $\mathfrak{g}_{0}$-module should be completely reducible. 
Even for the needs of the First Theorem of Invariant Theory this is
too strong a requirement, cf.  examples with complete irreducibility
in \cite{ser1, ser2} with our last example, in which the complete
reducibility of $\mathfrak{pe}(n)$ is violated.  Investigation of the
requiremets on $\Gamma$ and $\Gamma'$ needed for the First Theorem of
Invariant Theory will be given elsewhere.

{\bf 2.5.3.  Bernstein's square root of the Lefschetz decomposition}. 
Let $L$ be the space of a (complex) line bundle over a connected
symplectic manifold $(M^{2n}, \omega)$ with connection $\nabla$ such
that the curvature form of $\nabla$ is equal to $\hbar \omega$ for
some $\hbar \in{\mathbb C}$.  This $\hbar$ will be called a {\it
twist}; the space of tensor fields of type $\rho$ (here
$\rho:\mathfrak{sp}(2n)\longrightarrow \mathfrak{gl}(U)$ is a
representation which defines the space $\Gamma(M, U)$ of tensor fields
with values in $U$), and twist $\hbar$ will be denoted by
$T_\hbar(\rho)$.  Let us naturally extend the action of $X_+$, $X_-$
from the space $\Omega$ of differential forms on $M$ onto the space
$\Omega_\hbar$ of twisted differential forms using the isomorphism of
{\it spaces} $T_\hbar(\rho)\simeq T(\rho)\otimes \Gamma(L)$, where
$\Gamma(L)=\Omega^0_\hbar$ is the space of sections of the line bundle
$L$, i.e., the space of twisted functions.

Namely, set $X_+\mapsto X_+\otimes 1$, etc.  Let $D_+=d+\alpha$ be the
connection $\nabla$ itself and $D_-=[X_-, D_+]$.  On $\Omega_\hbar$,
introduce a superspace structure setting $p(\varphi\otimes s)
=\deg\varphi\pmod 2$, for $\varphi\in\Omega$, $s\in\Omega^0_\hbar$.

{\bf Theorem}.  (\cite{ber}) {\it On $\Omega_\hbar$, the operators
$D_+$ and $D_-$ generate an action of the Lie superalgebra
$\mathfrak{osp}(1|2)$ commuting with the action of the group $\hat G$
of $\nabla$-preserving automorphisms of the bundle $L$.}

Bernstein studied the $\hat G$-action, more exactly, the action of the
Lie algebra $\mathfrak{po}(2n|0)$ corresponding to $\hat G$; we are
interested in the part of this action only: in
$\mathfrak{sp}(2n)=\mathfrak{po}(2n|0)_0$-action.

In Example 2.5.1 the space $P^i$ consisted of differential forms with
constant coefficients.  Denote by ${\mathcal P}^i=P^i\otimes
S^{\bcdot}(V)$ the space of primitive forms with polynomial
coefficients.  The elements of the space $\sqrt{{\mathcal P}}^i_\hbar
=\ker D_-\cap {\mathcal P}^i_\hbar$ will be called $\nabla$-{\it
primitive forms} of degree $i$ (and twist $\hbar$).

Bernstein showed that $\sqrt{{\mathcal P}}^i_\hbar$ is an irreducible
$\mathfrak{g}=\mathfrak{po}(2n|0)$-module.  It could be that over
subalgebra $\mathfrak{g}_{0}$ the module $\sqrt{{\mathcal P}}^i_\hbar$
becomes reducible but the general theorem of Howe (which is true for
$\mathfrak{osp}(1|2n)$) states that this is not the case, it remains
irreducible.  Shapovalov and Shmelev literally generalized Bernstein's
result for $(2n|m)$-dimensional supermanifolds, see review
\cite{lei3}.  In particular, Shapovalov, who considered $n=0$, ``took
a square root of Laplacian and the metric''.

{\bf 2.5.4}.  Inspired by Bernstein's construction, let us similarly
define a ``square root'' of the hyper-K\"ahler structure.  Namely, on
a hyper-K\"ahlerean manifold $(M, \omega_1, \omega_2)$ consider a line
bundle $L$ with two connections: $\nabla_1$ and $\nabla_2$, whose
curvature forms are equal to $\hbar_1 \omega_1$ and $\hbar_2 \omega_2$
for some $\hbar_1, \hbar_2 \in{\mathbb C}$.  The pair $\hbar=(\hbar_1,
\hbar_2)$ will be called a {\it twist}; the space of tensor fields of
type $\rho$ and twist $\hbar$ will be denoted by $T_\hbar(\rho)$. 
Verbitsky \cite{ver} defined the action of $\mathfrak{sp}(4)$ in the
space $\Omega$ of differential forms on $M$.  Let us naturally extend
the action of the generators $X^{\pm}_j$ for $j= 1, 2$ of of
$\mathfrak{sp}(4)$ from $\Omega$ onto the space $\Omega_\hbar$ of
twisted differential forms using the isomorphism $T_\hbar(\rho)\simeq
T(\rho)\otimes \Gamma(L)$, where $\Gamma(L)=\Omega^0_\hbar$ is the
space of sections of the line bundle $L$; here $X^{+}_j$ is the
operator of multiplication by $\omega_j$ and $X^{-}_j$ is the operator
of convolution with the dual bivector.

Define the space of primitive $i$-forms (with constant coefficients)
on the hyper-K\"ahlerean manifold $(M, \omega_1, \omega_2)$ by setting
$$
P^i =\ker X^-_1\cap \ker X^-_2\cap\Omega ^i.\eqno{(HK)}
$$
According to the 
general theorem \cite{how2} this space is an irreducible 
$\mathfrak{sp}(2n; {\mathbb H})$-module.

Set $D^-_i=[X^-_i, D^+_i]$. The promised square root of the 
decomposition (HK) is the space
$$
{\mathcal P}^i_\hbar =\ker D^-_1\cap \ker D^-_2\cap\Omega 
^i_\hbar. \eqno{(\sqrt{{\rm HK}})}
$$
The operators $D^{\pm}_i$, where $D^+_i=\nabla_{i}$, generate
$\mathfrak{osp}(1|4)$.

{\bf 2.6. Further examples of dual pairs}.  The following subalgebras
$\mathfrak{g}_{1}(V_{1})\oplus \mathfrak{g} _{2}(V_{2})$ are maximal
in $\mathfrak{g}(V_{1}\otimes V_{2})$, hence, are dual pairs:
$$
\begin{array}{|c|c|c|}
\hline
\mathfrak{g}_{1} & \mathfrak{g}_{2} & \mathfrak{g} \\
\hline
\mathfrak{osp}(n_{1}|2m_{1}) & \mathfrak{osp} (n_{2}|2m_{2}) & 
\mathfrak{osp} (n_{1}n_{2}+4m_{1}m_{2}|2n_{1}m_{2}+2n_{2}m_{1}) \\
\mathfrak{o}(n) & \mathfrak{osp} (n_{2}|2m_{2}) & \mathfrak{osp} (nn_{2}
|2nm_{2}), n\neq 2, 4 \\
\mathfrak{sp}(2n) & \mathfrak{osp} (n_{2}|2m_{2}) & 
\mathfrak{osp} (2mn_{2}|4nm_{2}) \\
\mathfrak{pe} (n_{1}) & \mathfrak{pe}(n_{2}) & 
\mathfrak{osp} (2n_{1}n_{2}|2n_{1}n_{2}), n_{1}, n_{2}>2 \\
\hline
\mathfrak{osp} (n_{1}|2m_{1}) & \mathfrak{pe} (n_{2}) & 
\mathfrak{pe} (n_{1}n_{2}+2m_{1}n_{2})\; {\rm  if}\; n_{1}\neq 2m_{1} \\
& & \mathfrak{spe} (n_{1}n_{2}+2m_{1}n_{2})\; {\rm  if}\; n_{1}=2m_{1} \\
\mathfrak{o} (n) & \mathfrak{pe} (m) & \mathfrak{pe} (nm) \\
\mathfrak{sp} (2n) & \mathfrak{pe} (m) & \mathfrak{pe} (2nm) \\
\hline
\end{array}
$$ 
In particular, on the superspace of polyvector fields, there is a
natural $\mathfrak{pe}(n)$-module structure, and $\mathfrak{pe}(1)$,
its dual partner in $\mathfrak{osp} (2n|2n)$, is spanned by the
divergence operator $\Delta$ (``odd Laplacian''), called the {\it BRST
operator}\index{BRST operator}\index{operator, BRST} (\cite{bat/tyu}),
the even operator of $\mathfrak{pe}(1)$ being
$\deg_{x}-\deg_{\theta}$, where $\theta_i=\pi(\frac{\partial}{\partial
x_{i}})$, $\pi$ being the shift of parity operator.

For further examples of maximal subalgebras in $\mathfrak{gl}$ and
$\mathfrak{q}$ see \cite{Shc1}.  These subalgebras give rise to other
new examples of Howe dual pairs.  For the decomposition of the tensor
algebra corresponding to some of these examples see \cite{ser1, ser2},
some of the latter are further elucidated in \cite{che/wan}.  Some
further examples of Howe's duality, considered in a detailed version
of our lectures, are: (1) over reals; (2) dual pairs in simple
subalgebras of $\mathfrak{po}(2n|m)$ distinct from
$\mathfrak{osp}(m|2n)$; in particular, (3) embeddings into
$\mathfrak{po}(2n|m; r)$, the nonstandard regradings of the Poisson
superalgebra, cf.  \cite{Shc2}; (4) a ``projective'' version of the
Howe duality associated with embeddings into the Lie superalgebra of
Hamiltonian vector fields, the quotient of the Poisson superalgebra,
in particular, the exceptional cases in dimension $(2|2)$, cf. 
\cite{lei/shc}.  It is also interesting to consider the prime
characteristic and an ``odd'' Howe's duality obtained from
quantization of the antibracket (the main objective of
\cite{del/lei}), to say nothing of $q$-quantized versions of the
above.

\section*{\S 3. Generalities on spinor and spinor-like representations}

{\bf 3.1.  The spinor and oscillator representations of Lie algebras}. 
The importance of the spinor representation became clear very early. 
One of the reasons is the following.  As is known from any textbook on
representation theory, the fundamental representations $R(\varphi_1)=
W$, $R(\varphi_2)=\Lambda ^2(W)$,\dots , $R(\varphi_{n-1})=\Lambda
^{n-1}(W)$ of $\mathfrak{sl}(W)$, where $\dim W=n$ and $\varphi_{i}$
is the highest weight of $\Lambda ^{i}(W)$, are irreducible.  Any
finite dimensional irreducible $\mathfrak{sl}(n)$-module $L^\lambda$
is completely determined by its highest weight $\lambda=
\sum\lambda_i\varphi_i$ with $\lambda_i\in{\mathbb Z}_+$.  The module
$L^\lambda$ can be realized as a submodule (or quotient) of
$\mathop{\otimes}\limits
\left(R(\varphi_i)^{\otimes\lambda_{i}}\right)$.

Similarly, every irreducible $\mathfrak{gl}(n)$-module $L^\lambda$,
where $\lambda= (\lambda_1, \dots \lambda_{n-1}; c)$ and $c$ is the
eigenvalue of the unit matrix, is realized in the space of tensors,
perhaps, twisted with the help of $c$-densities, namely in the space
$\mathop{\otimes}\limits_i
\left(R(\varphi_i)^{\otimes\lambda_{i}}\right )\otimes {\rm tr}^c$,
where ${\rm tr}^c$ is the Lie algebraic version of the $c$th power of
the determinant, i.e., infinitesimally, trace, given for any
$c\in{\mathbb C}$ by the formula $X\mapsto c\cdot {\rm tr} (X)$ for
any matrix $X\in \mathfrak{gl}(W)$.  Thus, all the irreducible finite
dimensional representations of $\mathfrak{sl}(W)$ are naturally
realized in the space of tensors, i.e., in the subspaces or quotient
spaces of the space $T^p_{q}=\underbrace{W\otimes \dots\otimes W}_{p}
\underbrace{\otimes W^{*}\otimes\dots \otimes W^{*}}_{q}$, where $W$
is the space of the identity representation.  For
$\mathfrak{gl}(W)$, we have to consider the space $T^p_{q}\otimes {\rm
tr}^c$.

For $\mathfrak{sp}(W)$, the construction is similar, except the
fundamental module $R(\varphi_i)$ is now a {\it part} of the module
$\Lambda ^i({\rm id})$ consisting of the {\it primitive forms}.

For $\mathfrak{o}(W)$, the situation is totally different: not all
fundamental representations can be realised as (parts of) the modules
$\Lambda ^i({\rm id})$.  The exceptional one (or two, for
$\mathfrak{o}(2n)$) of them is called the {\it spinor
representation}\index{spinor representation}\index{representation,
spinor}; for $\mathfrak{o}(W)$, where $\dim W=2n$, it is realized in
the Grassmann algebra $E^{\bcdot}(V)$ of a ``half'' of $W$, where
$W=V\oplus V^{*}$ is a decomposition into the direct sum of subspaces
isotropic with respect to the form preserved by $\mathfrak{o}(W)$. 
For $\dim W=2n+1$, it is realized in the Grassmann algebra
$E^{\bcdot}(V\oplus W_{0})$, where $W=V\oplus V^{*}\oplus W_{0}$ and
$W_{0}$ is the 1-dimensional space on which the orthogonal form is
nondegenerate.

The quantization of the harmonic oscillator leads to an infinite
dimensional analog of the spinor representation which after Howe we
call {\it oscillator}\index{oscillator
representation}\index{representation, oscillator} representation of
$\mathfrak{sp}(W)$.  It is realized in $S^{\bcdot}(V)$, where as
above, $V$ is a maximal isotropic subspace of $W$ (with respect to the
skew form preserved by $\mathfrak{sp}(W)$).  The remarkable likeness
of the spinor and oscillator representations was underlined in a
theory of {\it dual Howe's pairs},\index{dual Howe's pairs}
\cite{How}.

The importance of spinor-oscillator representations is different for
distinct classes of Lie algebras and their representations.  In the
description of irreducible finite dimensional representations of
classical matrix Lie algebras $\mathfrak{gl}(n)$, $\mathfrak{sl}(n)$
and $\mathfrak{sp}(2n)$ we can do without either spinor or oscillator
representations.  We can not do without spinor representation for
$\mathfrak{o}(n)$, but a pessimist might say that spinor
representation constitutes only $\frac1n$th of the building bricks. 
Our, optimistic, point of view identifies the spinor representations
as one of the two possible types of the building bricks.

For the Witt algebra $\mathfrak{witt}$ and its central extension, the
Virasoro algebra $\mathfrak{vir}$, {\it every} irreducible highest
weight module is realized as a quotient of a spinor or, equivalently,
oscillator representation, see \cite{Fei/Fuc}, \cite{Fei/Fre}.  This
miraculous equivalence is known in physics under the name of {\it
bose-fermi correspondence}\index{bose-fermi
correspondence}\index{correspondence, bose-fermi}, see
\cite{Gre/Sch/Wit}, \cite{Kac}.

For the list of generalizations of $\mathfrak{witt}$ and
$\mathfrak{vir}$, i.e., simple (or close to simple) stringy Lie
superalgebras or Lie superalgebras of vector fields on $N$-extended
supercircles, often called by an unfortunate (as explained in  
\cite{Gro/Lei/Shc}) name ``superconformal
algebras'', see \cite{Gro/Lei/Shc}.  The importance of spinor-oscillator
representations diminishes as $N$ grows, but for the most interesting
 --- {\it distinguished} (\cite{Gro/Lei/Shc}) --- stringy superalgebras it is high, cf.
\cite{Fei/Sem/Sir/Tip}, \cite{Pol}.

{\bf 3.2.  Semi-infinite cohomology}.  An example of applications of
spinor-oscillator representations: semi-infinite (or BRST) cohomology
of Lie superalgebras.  These cohomology were introduced by Feigin
first for Lie algebras (\cite{Fei/Lei}); then he extended the
definition to Lie superalgebras via another construction, equivalent
to the first one for Lie algebras (\cite{Fei}).  For an elucidation of
Feigin's construction see \cite{Fre/Gar/Zuc}, \cite{Kos/Ste} and
\cite{Vor}.  Feigin rewrote in mathematical terms and generalized the
constructions physicists used to determine the {\it critical
dimensions} of string theories, i.e., the dimensions in which the
quantization of the superstring is possible, see \cite{Mar},
\cite{Gre/Sch/Wit}.  These critical dimensions are the values of the
central element (central charges) on the spinor-oscillator
representation constructed from the adjoint representation; to this
day not for every central element of all distinguished simple stringy
superalgebras their values are computed on every spinor-oscillator
representation, not even on the ones constructed from the adjoint
representations.

\section*{\S 4. The spinor-oscillator representations and Lie
superalgebras}

{\bf 4.1.  Spinor (Clifford--Weil--wedge--\dots) and oscillator
representations}.  As we saw in \cite{lei/shc},
$\mathfrak{po}(2n|m)_0\cong \mathfrak{osp}(m|2n)$, the superspace of
elements of degree 0 in the standard ${\mathbb Z}$-grading of
$\mathfrak{po}(2n|m)$ or, which is the same, the superspace of
quadratic elements in the representation by generating functions.  At
our first lecture we defined the {\it
spinor-oscillator representation}\index{spinor-oscillator
representation}\index{representation, spinor-oscillator} as the
through map (here $k=[\frac{m}{2}]$ and ${\mathcal Q}$ is the
quantization)
$$
\mathfrak{g}\longrightarrow \mathfrak{po}(2n|m)\stackrel{{\mathcal
Q}}{\longrightarrow} \left\{\begin{array}{l} 
\mathfrak{diff}(n|k)\quad {\rm if }\quad m=2k \\
\mathfrak{qdiff}(n|k)\quad {\rm if } \quad m=2k-1,\end{array}\right.
$$
where ${\rm Im}(\mathfrak{g})\subset
\mathfrak{po}(2n|m)_{0}=\mathfrak{osp}(m|2n)$.  Actually, such
requirement is too restrictive, we only need that the image of
$\mathfrak{g}$ under embedding into $\mathfrak{po}(2n|m)$ remains
rigid under quantization.  So various simple subalgebras of
$\mathfrak{po}(2n|m)$ will do as ambients of $\mathfrak{g}$.

This spinor-oscillator representation is called the {\it spinor
representation} of $\mathfrak{g}$ \index{spinor representation}
\index{representation, spinor} if $n=0$, or the {\it oscillator
representation}\index{oscillator representation}\index{representation,
oscillator} if $m=0$.  We will denote this representation ${\rm Spin}
(V)$ and set ${\rm Osc}(V)={\rm Spin} (\Pi(V))$, where $V$ is the
standard representation of $\mathfrak{osp}(m|2n)$.  In other words, if
${\rm Spin} (V)$ is a representation of $\mathfrak{osp}(m|2n)$, then
${\rm Osc}(V)$ is a representation of $\mathfrak{osp}(2n|m)$, so ${\rm
Osc}(V)$ only exists for $m$ even.

If $V$ is a $\mathfrak{g}$-module without any bilinear form, but we
still want to construct a spinor-oscillator representation of
$\mathfrak{g}$, consider the module $W=V\oplus V^*$ (where in the
infinite dimensional case we replace $V^*$ with the {\it restricted}
dual of $V$; roughly speaking, if $V={\mathbb C}[x]$, then
$V^{*}={\mathbb C}[[\frac{\partial}{\partial x}]]$, whereas the
restricted dual is ${\mathbb C}[\frac{\partial}{\partial x}]$) endowed
with the form (for $v_1, w_1\in V$, $v_2, w_2\in V^*$) symmetric for
the plus sign and skew-symmetric otherwise:
$$
B((v_1, v_2), (w_1, w_2))= v_2(w_1)\pm
(-1)^{p(v_{1})p(w_{2})}w_2(v_1).
$$
Now, in $W$, select a maximal isotropic subspace $U$ (not necessarily $V$ 
or $V^*$) and realize the spinor-oscillator representation of
$\mathfrak{g}$ in the exterior algebra of $U$.

Observe that the classical descriptions of spinor representations
differ from ours, see, e.g., \cite{Got/Gro}, where the embedding of
$\mathfrak{g}$ (in their case $\mathfrak{g}=\mathfrak{o}(n)$) into the
quantized algebra (namely into ${\mathcal Q}(\mathfrak{po}(0|n-1))$)
is considered, not into $\mathfrak{po}(0|m)$.  The existence of this
embedding is not so easy to see unless told, whereas our constructions
are manifest and bring about the same result.

To illustrate our definitions and constructions, we realize the
orthogonal Lie algebra $\mathfrak{o}(n)$ as the subalgebra in the Lie
superalgebra $\mathfrak{po}(0|n)$.

Case $\mathfrak{o}(2k)$. Basis:
$$
\renewcommand{\arraystretch}{1.4}
\begin{array}{cccc}
X^+_1=\xi_2\eta_1,& \dots ,& X^+_{k-1}=\xi_{k}\eta_{k-1},&
X^+_{k}=\eta_{k}\eta_{k-1};\\
X^-_1=\xi_1\eta_2,& \dots ,& X^-_{k-1}=\xi_{k-1}\eta_{k},&
X^-_{k}=\xi_{k-1}\xi_{k};\\
H_1=\xi_1\eta_1-\xi_2\eta_2,& \dots ,&
H_{k-1}=\xi_{k-1}\eta_{k-1}-\xi_{k}\eta_{k},&
H_{k}=\xi_{k-1}\eta_{k-1}+\xi_{k}\eta_{k}.\end{array}
$$
For $ R(\varphi_k)$ take the subspacespace functions ${\mathbb
C}[\xi]_{{\rm ev}}$ which contains the constants ${\mathbb C}\cdot \hat 1$, where
$\hat 1$ is just the constant function $1$; clearly, $\hat 1$ is the
vacuum vector.

Quantization (see above) sends: $\xi_i$ into
$\hat\xi_i$, and $\eta_i$ into $\hbar\frac{\partial}{\partial\xi_i}$,
so $X^{\pm}_i\hat 1=0$ for $i<k$, hence, $H_i\hat 1=[X^+_i, X^-_i]\hat
1=0$ for $i<k$.  Contrariwise,
$$
 H_k\hat 1=[X^+_k, X^-_k]\hat 1=[\partial_k\partial_{k-1},
 \hat\xi_{k-1} \hat\xi_k]\hat 1= \partial_k(-
 \hat\xi_{k-1}\partial_{k-1}+1) \hat\xi_k\hat 1=\hat 1.
$$
So we see that the spinor representation is indeed a fundamental one. 

Case $\mathfrak{o}(2k+1)$. Basis:
$$
\renewcommand{\arraystretch}{1.4}
\begin{array}{cccc}
X^+_1=\xi_2\eta_1,& \dots ,& X^+_{k-1}=\xi_{k}\eta_{k-1},&
X^+_{k}=\sqrt{2}\eta_{k}\theta;\\
X^-_1=\xi_1\eta_2,& \dots ,& X^-_{k-1}=\xi_{k-1}\eta_{k},&
X^-_{k}=\sqrt{2}\theta\xi_{k};\\
H_1=\xi_1\eta_1-\xi_2\eta_2,& \dots ,&
H_{k-1}=\xi_{k-1}\eta_{k-1}-\xi_{k}\eta_{k},&
H_{k}=2\xi_{k}\eta_{k}.\end{array}
$$
For $R(\varphi_k)$ consider the space of even functions ${\mathbb
C}[\xi_1, \dots , \xi_k, \theta]_{{\rm ev}}$ and realize
$\mathfrak{o}(2k+1)$ so that $\xi_i\mapsto \hat\xi_i$, $\eta_i\mapsto
\hbar\frac{\partial}{\partial\xi_i}$, $\theta\mapsto
\hbar( \hat\theta+\frac{\partial}{\partial\theta})$.  As above for
$\mathfrak{o}(2k)$, set $\hbar=1$.

Then, as above, $H_iv=[X^+_i, X^-_i]\hat 1=0$ for $i<k$, whereas
$$
H_k\hat 1=[X^+_k, X^-_k]\hat 1=\frac 22\left(\partial_k(
\hat\theta+\frac{\partial}{\partial\theta})^2 \hat\xi_k+ \hat\xi_k(
\hat\theta+\frac{\partial}{\partial\theta})^2\partial_k\right)\hat
1=\hat 1.
$$ 
So $\hat 1$ is indeed the highest weight vector of the $k$th
fundamental representation.  

{\bf 4.2.  Stringy superalgebras. Case $\mathfrak{vir}$}. For the
basis of $\mathfrak{vir}$ take $e_{i}=t^{i+1}\frac{d}{dt}$,
$i\in{\mathbb Z}$, and the central element $z$; let the bracket be 
$$
[e_{i}, e_{j}]=(j-i)e_{i+j}-\frac{1}{12}\delta_{ij}(i^3-i)z. \eqno{(*)}
$$
We advise the reader to refresh definitions of stringy superalgebras
and various modules over them, see \cite{Gro/Lei/Shc}, where we also
try to convince physicists not to use the term ``superconformal
algebra'' (except, perhaps, for $\mathfrak{k}^L(1|1)$ and
$\mathfrak{k}^M(1|1)$).  In particular, recall that ${{\mathcal
F}}_{\lambda, \mu}= {\rm Span}(\varphi_{i}=t^{\mu+i}(dt)^{\lambda}\mid
i\in{\mathbb Z})$.

{\bf Statement}.  {\it The only instances when ${{\mathcal
F}}_{\lambda, \mu}$ possesses an invariant {\em symmetric}
nondegenerate bilinear form are the space of half-densities,
$\sqrt{{\rm Vol}}={\mathcal F}_{1/2, 0}$, and its twisted version,
${\mathcal F}_{1/2, 1/2}$ and in both cases the form is:
$$
(f\sqrt{dt}, g\sqrt{dt})=\int fg\cdot dt;
$$
the only instances when ${{\mathcal F}}_{\lambda, \mu}$ possesses an
invariant {\em skew-symmetric} forms are the quotient space of
functions modulo constants, $d{\mathcal F}={\mathcal F}_{0,
0}/{\mathbb C}\cdot 1$, and $\frac{1}{2}$-twisted functions,
$\sqrt{t}{\mathcal F}={\mathcal F}_{0, 1/2}$ and in both cases the
form is:
$$
(f, g)=\int f\cdot dg.
$$
}

Let $\partial_i=\frac{\partial}{\partial \varphi_i}$ (where
$\varphi_{i}=t^{\mu+i}(dt)^{\lambda}$).  Let ${\rm osc} (\sqrt{{\rm
Vol} })$ be the $\mathfrak{vir}$-submodule of the {\it exterior}
algebra on $\varphi_i$ for $i<0$ 
containing the constant $\hat 1$.  Since the generators $e_i$ of
$\mathfrak{vir}$ acts on ${\mathcal F}_{\lambda, \mu}$ as (sums over 
$i\in {\mathbb Z}$)
$$ 
\renewcommand{\arraystretch}{1.4}
\begin{array}{rl}
e_1=&\sum(\mu+i+2\lambda)\varphi_{i+1}\partial_i=\sum i\varphi_{i+1}\partial_i,\\
e_{-1}=&\sum(\mu+i+1)\varphi_{i}\partial_{i+1}=\sum
(i+1)\varphi_{i}\partial_{i+1};\\
e_2=&\sum(\mu+i-\lambda)\varphi_{i+1}\partial_i=\sum i\varphi_{i+1}\partial_i,\\
e_{-2}=&\sum(\mu+i+3\lambda)\varphi_{i}\partial_{i+1}=\sum
(i+1)\varphi_{i}\partial_{i+1},
\end{array}
$$
and representing $e_0$ and $z$ as brackets of $e_{\pm 1}$ and $e_{\pm
2}$ from $(*)$ we immediately deduce that the highest weights $(c, h)$
of ${\rm osc} (\sqrt{{\rm Vol} })$ is $(-\frac{1}{ 3}, 0)$.

For the spinor representations ${\rm spin} (\sqrt{t}{{\mathcal F}})$
and ${\rm spin} (d{\mathcal F})$ (realized on the {\it symmetric}
algebra of $\varphi_i$ for $i<0$) we similarly obtain that the highest
weights $(c, h)$ are $(\frac{1}{6}, \frac{1}{2})$ for ${\rm spin}
(\sqrt{t}{\mathcal F})$ and $(-\frac{1}{6}, 0)$ for ${\rm spin}
(d{\mathcal F})$.

Observe that the representations ${\rm spin} (\sqrt{t}{\mathcal F})$,
${\rm spin} (d{\mathcal F})$ and ${\rm osc} (\sqrt{{\rm Vol} })$ are
constructed on a half of the generators used to construct ${\rm Spin}
({\mathcal F}_{\lambda; \mu})$.

{\bf 4.3.  The highest weights of the spinor representations of
$\mathfrak{k}^L(1|n)$ and $\mathfrak{k}^M(1|n)$}.  In the following
theorem we give the coordinates $(c, h; H_1, \dots )$ of the highest
weight of the spinor representations ${\rm Spin} ({\mathcal
F}_{\lambda; \mu})$ of the contact superalgebra $\mathfrak{k}^L(1|n)$
with respect to $z$ (the central element), $K_t$, and, after
semicolon, on the elements of Cartan subalgebra, respectively.  For
$\mathfrak{k}^M(1|n)$ we write $\tilde h;\tilde H_i$.  (Observe that
for $n>4$ the Cartan subalgebra has more generators than just $H_1=
K_{\xi_{1}\eta_{1}}$, \dots , $H_k=K_{\xi_{k}\eta_{k}}$ which generate
the Cartan subalgebra of $\mathfrak{k}(1|2k)$, the algebra of contact
vector fields with polynomial coefficients.)

\noindent \begin{tabular}{|c|c|c|c|c|}
\hline
$n$&$ 0 $&$1  $&$2  $&$\geq 3$\\
\hline
$c$&$12\lambda^2-12\lambda+2$&$-12\lambda+3$&$6$&$0$\\
\hline
$h$&$(\mu+2\lambda)(\mu+1) $&$\mu+2\lambda $&$ 2\mu+2\lambda +\nu$&$
2^{n-1}(\mu+\lambda)+2^{n-3}$\\
\hline
$\tilde h$&--&$2\mu+3\lambda-\frac{1}{4}$&$2\mu+2\lambda-\frac{1}{2}$&
$2^{n-1}(\mu+\lambda)$\\
\hline
\end{tabular}

{\bf Theorem}.  {\sl Let $(c, h; H_1, \dots )$ be the highest weight
of the spinor representation ${\rm Spin} ({\mathcal F}_{\lambda;
\mu})$ of $\mathfrak{k}^L(1|n)$.  The highest weight of the oscillator
representation ${\rm Osc} ({\mathcal F}_{\lambda; \mu})={\rm Spin}
(\Pi({{\mathcal F}}_{\lambda; \mu}))$ is $(-c, h; H_1, \dots)$ and
similarly for $\mathfrak{k}^M(1|n)$.

For $n\neq 2$, all the coordinates of the highest weight other than
$c$, $h$ vanish.  For $n=2$ the value of $H$ on the highest weight
vector from ${\rm Spin} ({{\mathcal F}}_{\lambda, \nu; \mu})$ is equal
to $\nu$.

The values of $c$ and $h$ (or $\tilde h$) on modules ${\rm Spin}
({{\mathcal F}}_{\lambda; \mu})$ are given in the above table.}

\noindent Up to rescaling, these results are known for small $n$, see \cite{Ken/Rig},
\cite{Kac/Wak} and refs. therein.

{\bf Remark}.  For the contact superalgebras $\mathfrak{g}$ on the
$1|n$-dimensional supercircle our choice of $\mathfrak{g}$-modules
$V={\mathcal F}_{\lambda; \mu}$ from which we constructed ${\rm Spin}
(V\oplus V^*)$ is natural for small $n$: there are no other modules! 
For larger $n$ it is only justified if we are interested in
semi-infinite cohomology of $\mathfrak{g}$ and not in representation
theory {\it per se}.  For the superalgebras $\mathfrak{g}$ of series
$\mathfrak{vect}$ and $\mathfrak{svect}$ the adjoint module
$\mathfrak{g}$ is of the form ${\mathcal T}({\rm id}^*)$, i.e, it is
either coinduced from multidimensional representation
($\mathfrak{vect}$), or is a submodule of such a coinduced module
($\mathfrak{svect}$). Spinor-oscillator representations of this type
were not studied yet, cf.  sec.  5.

{\bf 4.4.  Other spinor representations}. 1) Among various Lie
superalgebras for which it is interesting to study spinor-oscillator
representations, the simple (or close to them) maximal
subsuperalgebras of $\mathfrak{po}$ are most interesting.  The list of
such maximal subalgebras is being completed; various maximal
subalgebras listed in \cite{shc} distinct from the sums of mutual
centralizers also provide with spinor representations. 

As an interesting example consider A.~Sergeev's Lie superalgebra
$\mathfrak{as}$, the nontrivial central extension of the Lie
superalgebra $\mathfrak{spe}(4)$ preserving the odd bilinear form and
the volume on the $(4|4)$-dimensional superspace, see \cite{Shc1,
Shc2}.  Namely, consider $\mathfrak{po}(0|6)$, the Lie superalgebra
whose superspace is the Grassmann superalgebra $\Lambda(\xi, \eta)$
generated by $\xi_1, \xi_2, \xi_3, \eta _1, \eta _2, \eta_3$ and the
bracket is the Poisson bracket.  Recall also that the quotient of 
$\mathfrak{po}(0|6)$ modulo center is $\mathfrak{h}(0|6)={\rm
Span} (H_f \mid f\in\Lambda (\xi, \eta))$, where
$$
H_f=(-1)^{p(f)}\sum(\frac{\partial f}{\partial
\xi_j} \frac{\partial}{\partial\eta_j}+ \frac{\partial
f}{\partial\eta_j} \frac{\partial}{\partial\xi_j}).
$$
Now, observe that $\mathfrak{spe}(4)$ can be embedded into
$\mathfrak{h}(0|6)$.  Indeed, setting $\deg \xi_i=\deg \eta _i=1$ for
all $i$ we introduce a ${\mathbb Z}$-grading on $\Lambda(\xi, \eta)$
which, in turn, induces a ${\mathbb Z}$-grading on $\mathfrak{h}(0|6)$
of the form $\mathfrak{h}(0|6)=\mathop{\oplus}\limits_{i\geq
-1}\mathfrak{h}(0|6)_i$.  Since
$\mathfrak{sl}(4)\cong\mathfrak{o}(6)$, we can identify
$\mathfrak{spe}(4)_0$ with $\mathfrak{h}(0|6)_0$.

It is not difficult to see that the elements of degree $-1$ in the
standard gradings of $\mathfrak{spe}(4)$ and $\mathfrak{h}(0|6)$
constitute isomorphic $\mathfrak{sl}(4)\cong\mathfrak{o}(6)$-modules. 
It is subject to a direct verification that it is really possible to
embed $\mathfrak{spe}(4)_1$ into $\mathfrak{h}(0|6)_1$.

A.~Sergeev's extension $\mathfrak{as}$ is the result of the
restriction onto $\mathfrak{spe}(4)\subset\mathfrak{h}(0|6)$ of the
cocycle that turns $\mathfrak{h}(0|6)$ into $\mathfrak{po}(0|6)$.  The
quantization (with parameter $\lambda$) deforms $\mathfrak{po}(0|6)$
into $\mathfrak{gl}(\Lambda(\xi))$; the through maps $T_\lambda:
\mathfrak{as}\longrightarrow\mathfrak{po}(0|6)\longrightarrow\mathfrak{gl}(\Lambda
(\xi))$ are representations of $\mathfrak{as}$ in the
$4|4$-dimensional modules ${\rm Spin}_\lambda$.  The explicit form of
$T_\lambda$ is as follows:
$$
T_\lambda: \left (\begin{array}{ll}a & b \\ 
c & -a^t \end{array} \right) +d\cdot z\mapsto  
\left (\begin{array}{ll}
a & b-\lambda \tilde c \\ 
c & -a^t \end{array} \right)+\lambda d\cdot 1_{4|4}, 
$$
where $1_{4|4}$ is the unit matrix and $\, \tilde{}\, $ is extended
via linearity from matrices $c_{ij}=E_{ij}-E_{ji}$ on which $\tilde
c_{ij}=c_{kl}$ for any even permutation $(1234)\mapsto(ijkl)$. 
Clearly, $T_\lambda$ is an irreducible representation for any
$\lambda$ and $T_\lambda\not\simeq T_\mu$ for $\lambda\neq\mu$.

2) Maximal subalgebras (for further examples see \cite{shc}) and a
conjecture.  Let $V_1$ be a linear superspace of dimension $(r|s)$;
let $\Lambda (n)$ be the Grassmann superalgebra with $n$ odd
generators $\xi_1$, \dots , $\xi_n$ and
$\mathfrak{vect}(0|n)=\mathfrak{der}\Lambda (n)$ the Lie superalgebra
of vector fields on the $(0|n)$-dimensional supermanifold.

Let $\mathfrak{g}=\mathfrak{gl}(V_1)\otimes \Lambda
(n)\subplus\mathfrak{vect}(0|n)$ be the semidirect sum (the ideal at
the closed part of $\subplus$) with the natural action of
$\mathfrak{vect}(0|n)$ on the ideal $\mathfrak{gl}(V_1)\otimes \Lambda
(n)$.  The Lie superalgebra $\mathfrak{g}$ has a natural faithful
representation $\rho$ in the space $V=V_1\otimes \Lambda (n)$ defined
by the formulas
$$
\begin{array}{ll} 
\rho(X\otimes \varphi)(v\otimes\psi)=(-1)^{p(\varphi)p(\psi)}Xv\otimes
\varphi\psi,\\
\rho(D)(v\otimes\psi)=-(-1)^{p(D)p(v)}v\otimes
D\psi\end{array}
$$
for any $X\in\mathfrak{gl}(V_1)$, $\varphi, \psi\in \Lambda (n)$,
$v\in V_1$, $D\in\mathfrak{vect}(0|n)$.  Let us identify the elements
from $\mathfrak{g}$ with their images under $\rho$, so we consider
$\mathfrak{g}$ embedded into $\mathfrak{gl}(V)$.

{\bf Theorem} (\cite{shc}) $1)$ {\sl The Lie superalgebra
$\mathfrak{gl}(V_1)\otimes \Lambda (n)\subplus \mathfrak{vect}(0|n)$
is maximal irreducible in $\mathfrak{sl}(V_1\otimes \Lambda (n))$
unless {\em a)} $\dim V_1=(1, 1)$ or {\em b)} $n=1$ and $\dim V_1=(1,
0)$ or $(0, 1)$ or $(r|s)$ for $r\neq s$.

$2)$ If $\dim V_1=(1, 1)$, then $\mathfrak{gl}(1|1)\cong \Lambda
(1)\subplus \mathfrak{vect}(0|1)$, so 
$$
\mathfrak{gl}(V_1)\otimes \Lambda (n)\subplus
\mathfrak{vect}(0|n)\subset \Lambda (n+1)\subplus
\mathfrak{vect}(0|n+1)
$$ 
and it is the bigger superalgebra which is maximal irreducible in
$\mathfrak{sl}(V)$.

$3)$ If $n=1$ and $\dim V_1= (r|s)$ for $r> s>0$, then
$\mathfrak{g}$ is maximal irreducible in $\mathfrak{gl}(V)$.}

{\bf Conjecture}.  Suppose $r+s=2^N$.  Then, $\dim V$ coincides with
$\dim \Lambda (W)$ for some space $W$.  We suspect that this
coincidence is not accidental but is occasioned by the spinor
representations of the maximal subalgebras described above.  The same
applies to $\mathfrak{q}(V_1)\otimes
\Lambda(n)\subplus\mathfrak{vect}(0|n)$, a maximal irreducible
subalgebra in $\mathfrak{q}(V_1\otimes \Lambda (n))$.

{\bf 4.5.  Selected problems}.  1) The spinor and oscillator
representations are realized in the symmetric (perhaps,
supersymmetric) algebra of the maximal isotropic (at least for
$\mathfrak{g}=\mathfrak{sp}(2k)$ and $\mathfrak{o}(2k)$) subspace $V$
of the identity $\mathfrak{g}$-module ${\rm id}=V\oplus V^*$.  But one
could have equally well started from another $\mathfrak{g}$-module. 
For an interesting study of spinor representations constructed from
$W\neq {\rm id}$, see \cite{Pan}.

To consider in a way similar to sec. 2 contact stringy
superalgebras $\mathfrak{g}=\mathfrak{k}^L(1|n)$ and
$\mathfrak{k}^M(1|n)$, as well as other stringy superalgebras from the
list \cite{Gro/Lei/Shc}, we have to replace ${{\mathcal F}}_{\lambda,
\mu}$ with modules ${{\mathcal T}}_{\mu}(W)$ of (twisted) tensor
fields on the supercircle and investigate how does the highest weight
of $\hat 1\in {\rm Osc} ({{\mathcal T}}_{\mu}(W))$ or $\hat 1\in {\rm
Spin} ({{\mathcal T}}_{\mu}(W))$ constructed from an arbitrary
irreducible $\mathfrak{co}(n)$-module $W=V\oplus V^*$ depend on the
highest weight of $W$.  (It seems that the new and absolutely
remarkable spinor-like representation Poletaeva recently constructed
\cite{Pol} is obtained in this way.)

To give the reader a feel of calculations, we consider here the
simplest nontrivial case $\mathfrak{o}(3)=\mathfrak{sl}(2)$.  The
results may (and will) be used in calculations of ${\rm Spin}
({{\mathcal T}}_{\mu}(W))$ for $\mathfrak{g}=\mathfrak{k}^L(1|n)$ and
$\mathfrak{k}^M(1|n)$ for $n=3, 4$.  As is known, for every
$N\in{\mathbb Z}_+$ there exists an irreducible $(N+1)$-dimensional
$\mathfrak{g}$-module with highest weight $N$.  This module possesses
a natural nondegenerate $\mathfrak{g}$-invariant bilinear form which
is skew-symmetric for $N=2k+1$ and symmetric for $N=2k$.  The
corresponding embeddings
$\mathfrak{g}\longrightarrow\mathfrak{o}(2k+1)$ and
$\mathfrak{g}\longrightarrow\mathfrak{sp}(2k)$ are called {\it
principal}, see \cite{Gro/Lei1} and references therein.  Explicitly,
the images of the Chevalley generators $X^{\pm}$ of $\mathfrak{sl}(2)$
are as follows: $X^-\mapsto\sum X^-_i$,
$$
\renewcommand{\arraystretch}{1.4}
X^+\mapsto \left\{\begin{array}{ll} N(N+1)X^+_N +
\sum\limits_{1\leq i\leq N-1} i(N+1-i)X^+_i&{\rm  for }\; N=2k+1 \\
N^2X^+_N + \sum\limits_{1\leq i\leq N-1} i(2N-i)X^+_i&{\rm  for }\;
N=2k .\end{array}\right.
$$

From the commutation relations between $X^{+}$ and $X^{-}$ we derive
that only $X^{\pm}_N$ give a nontrivial contribution to the highest
weight ${\rm HW}$ of the $\mathfrak{sl}(2)$-module ${\rm Spin}(L^N)$;
we have:
$$
\renewcommand{\arraystretch}{1.4}
{\rm HW}= \left\{\begin{array}{ll}N(N+1)&{\rm  if }\quad N=2k+1 \\
-\frac{1}{ 2} N^2&{\rm  if }\quad N=2k. \end{array}\right .
$$

2) Observe, that the notion of spinor-oscillator representation can be
broadened to embrace the subalgebras of the Lie superalgebra
$\mathfrak{h}$ of Hamiltonian vector fields and their images under
quantization; we call the through map the {\it projective
spinor-oscillator representation}.\index{projective spinor-oscillator
representation}\index{ representation spinor-oscillator, projective}
Since the Lie superalgebra $\mathfrak{h}$ has more deformations than
$\mathfrak{po}$ (\cite{lei/shc}), and since the sets of maximal simple
subalgebras of $\mathfrak{po}$ and $\mathfrak{h}$ are distinct, the
set of examples of projective spinor-oscillator representations
differs from that of spinor-oscillator representations.

\end{document}